\documentclass[11pt,oneside]{article}

\usepackage{amsmath}
\usepackage{amsthm}
\usepackage{amssymb}

\usepackage[all]{xy}
\SelectTips{cm}{12}
\CompileMatrices

\newtheorem{theorem}{Theorem}[section]
\newtheorem{proposition}[theorem]{Proposition}
\newtheorem{prop}[theorem]{Proposition}
\newtheorem{lemma}[theorem]{Lemma}
\newtheorem{corollary}[theorem]{Corollary}

\newtheorem{problem}[theorem]{Problem}

\theoremstyle{definition}
\newtheorem{defn}[theorem]{Definition}
\newtheorem{remark}[theorem]{Remark}
\newtheorem{convention}[theorem]{Convention}
\newtheorem{example}[theorem]{Example}

\def\bar{\overline}
\newcommand{\sdirect}{\rtimes}
\def\Bbb{\bf}

\newcommand\Z{\mbox{\Bbb Z}}
\newcommand\Q{\mbox{\Bbb Q}}
\newcommand\N{\mbox{\Bbb N}}

\DeclareMathOperator{\Fix}{Fix}

\renewcommand\to{\rightarrow}

\newcommand{\boundary}{\partial}

\DeclareMathOperator{\Aut}{Aut}
\DeclareMathOperator{\Vol}{Vol}

\newcommand{\A}{\mathcal{A}} 
\newcommand{\B}{\mathcal{B}} 
\newcommand{\inclusion}{\hookrightarrow}
\newcommand{\quotient}{\backslash\!\backslash}
\newcommand{\abs}[1]{\left| {#1} \right|}

\renewcommand{\emph}[1]{\textbf{#1}}
\def\endrem{}
	     
\usepackage{ifthen}

\newcommand{\showcomments}{yes}

\newsavebox{\commentbox}
{\unskip\ifthenelse{\equal{\showcomments}{yes}}%
{\footnotemark
    \begin{lrbox}{\commentbox}
    \begin{minipage}[t]{1.25in}\raggedright\sffamily\upshape\tiny
    \footnotemark[\arabic{footnote}]}
{\begin{lrbox}{\commentbox}}}%
{\ifthenelse{\equal{\showcomments}{yes}}%
{\end{minipage}\end{lrbox}\marginpar{\usebox{\commentbox}}}
{\end{lrbox}}}

\begin{document}

\title{Commensurability invariants for nonuniform tree lattices}

\author{Benson Farb\thanks{Supported in part by NSF grants DMS-9704640
and DMS-0244542.}\ \and
G.~Christopher Hruska\thanks{Supported in part by an NSF postdoctoral research
fellowship.}}




\date{\today}

\maketitle

\begin{abstract}
We study nonuniform lattices in the automorphism group $G$
of a locally finite simplicial tree~$X$.
In particular, we are interested in classifying
lattices up to commensurability in~$G$.
We introduce two new commensurability invariants:
\emph{quotient growth}, which measures the growth of the noncompact
quotient of the lattice; and \emph{stabilizer growth},
which measures the growth of the orders of finite stabilizers
in a fundamental domain as a function of distance from a fixed basepoint.
When $X$ is the biregular tree $X_{m,n}$,
we construct lattices realizing all triples of covolume,
quotient growth, and stabilizer growth satisfying some mild conditions.
In particular, for each positive real number~$\nu$ we construct
uncountably many noncommensurable lattices with covolume~$\nu$.
\end{abstract}



\section{Introduction}
\label{sec:Introduction}

Let $X$ be a locally finite simplicial tree.
Then the group $G=\Aut(X)$
of simplicial automorphisms of $X$ naturally has the structure of 
a locally compact topological group, with a neighborhood basis of the
identity consisting of automorphisms fixing larger and larger balls.  
A \emph{lattice} in $G$ is a
discrete subgroup of cofinite (Haar) measure $\mu$.
The study of these
``tree lattices'' generalizes the theory of lattices in 
rank one Lie groups over nonarchimedean local 
fields, and provides a remarkably rich theory (see the now-standard
reference \cite{BL}).
With the right normalization of the Haar
measure $\mu$, there is a combinatorial formula:
\[
   \Vol(\Gamma\quotient X):=\mu(\Gamma\backslash G)
      =\sum_{v\in V(A)}\frac{1}{\abs{\Gamma_v}}
\]
where the sum is taken over vertices $v$ in a fundamental domain
$A\subseteq X$ for the $\Gamma$--action,
and $\abs{\Gamma_v}$ is the order of
the $\Gamma$--stabilizer of~$v$.
Recall that a lattice $\Gamma<G$ is
\emph{uniform} if $\Gamma\backslash G$ is compact, and is 
\emph{nonuniform} otherwise. 

One of the basic problems about a locally compact topological group~$G$
is to
classify its lattices up to commensurability.
Recall that two lattices
$\Gamma_1,\Gamma_2 \le G$ are \emph{commensurable in $G$} if there exists
$g\in G$ so that $g\Gamma_1g^{-1}\cap \Gamma_2$ has finite index in
both $g\Gamma_1g^{-1}$ and $\Gamma_2$.
Since covolume is multiplicative
in index, two commensurable lattices have volumes which are
commensurable real numbers, \textit{i.e.}, which have a rational ratio.
Other
commensurability invariants are harder to come by unless $G$ is linear.

When $G=\Aut(X)$, Leighton proved in \cite{Leighton}
that all torsion-free uniform lattices are commensurable.
The torsion free hypothesis was later removed by Bass--Kulkarni in
\cite{BK}, establishing that there is at most
one commensurability class of uniform lattices in $G=\Aut(X)$.
The case of nonuniform
lattices, (which exist in abundance by a theorem of Carbone \cite{Ca}), 
is much more complicated.

In this paper we will concentrate on the case of the biregular tree
$X_{m,n}$ with degrees $m$ and~$n$.
Bass--Lubotzky \cite{BL} (for $m=n$)
and Rosenberg \cite{Ro} (for every $m\geq n\geq 3$) proved that for
every real number $r>0$, there exists a nonuniform lattice $\Gamma$ in
$G=\Aut(X_{m,n})$ with covolume~$r$.
This gives uncountably many
commensurability classes of nonuniform lattices in $G$, one for each
commensurability class of real numbers.

In Section~\ref{sec:Invariants} 
we introduce two commensurability invariants.  
The \emph{quotient growth type} (see \S\ref{subsec:QuotientGrowth})
is essentially the
``growth type'' of the quotient graph $\Lambda=\Gamma\backslash X$,
which is an equivalence class of functions 
measuring the growth of combinatorial balls in~$\Lambda$.
The \emph{stabilizer growth type} (see
\S\ref{subsec:StabilizerGrowth}) measures 
the growth of the order of the stabilizers of vertices in a fundamental
domain as a function of their distance to a fixed basepoint.
There is a refinement of this invariant for each prime $p$.

Our main result
(see Sections \ref{sec:construct} and~\ref{sec:construct2})
is the construction of lattices realizing uncountably many 
values of the invariants, 
further indicating the richness
of nonuniform lattices in $G$.
Informally, we prove:

\begin{theorem}[Main theorem, informal statement]
\label{theorem:informal}
Suppose $3 \le m \le n$, and let $G=\Aut(X_{m,n})$.
Let $r>0$ be any real number, and let $f,g\colon\N\to\N$ be any
functions satisfying the \textup{(}mild\textup{)} conditions given in 
Theorem \ref{theorem:newmain} below.  
Then there exists a nonuniform lattice $\Gamma<G$ with 
covolume $r$, quotient growth $f$,
and \textup{(}for $n>4$ composite\textup{)} stabilizer growth~$g$.
\end{theorem}

This result is stated precisely as Theorem \ref{theorem:newmain} below.
In particular we note the following.

\begin{corollary}
\label{corollary:main}
Let $m,n\geq 3$ and let $G=\Aut(X_{m,n})$.  Then for each real number
$r>0$, there exist uncountably many commensurability classes of 
nonuniform lattices in~$G$, each having covolume~$r$.
\end{corollary}

As a concrete example, given any $r>0$ and $1<\alpha <2$, we then construct
lattices $\Gamma_{r,\alpha}<G$ each with covolume $r$, and with the growth
type of the quotient tree proportional to $\alpha^k$; the number
$\alpha$ is a commensurability invariant.
Similarly, we can build
lattices whose quotients have ``intermediate growth,'' \textit{i.e.}, whose
$R$-balls have growth type $e^{R^\beta}$
for arbitrary $\beta \in (0,1)$.
In contrast note that, for
lattices in rank one Lie groups over nonarchimedean local fields, it
follows from theorems of Lubotzky
and Raghunathan (see \cite{Lu}) that such lattices always have quotients
with ``linear growth type.''

We would like to pose the following:

\begin{problem}
Classify nonuniform lattices in $\Aut(X_{m,n})$ up to commensurability.
\end{problem}

As a special case, it would be interesting to give a commensurability
classification for lattices all having a fixed quotient, for example a
ray.

\subsection{Acknowledgements}

The authors would like to thank the referee for several helpful comments
that we hope have improved the exposition of the paper.

\section{Preliminaries}
\label{sec:Preliminaries}

In this section we present some needed background material.

\subsection{Edge-indexed graphs and graphs of groups}
\label{EdgeIndexed}

In this subsection we briefly recall the basic tools in constructing tree
lattices.
We refer the reader to the book \cite{BL} by Bass--Lubotzky for a more
comprehensive introduction to this material.

A \emph{graph of groups} $\mathbf{A} = (A, \A)$ consists of a graph~$A$
with vertex set $V(A)$ and edge set $E(A)$, vertex groups $\A_v$
for each $v \in V(A)$, edge groups $\A_e = \A_{\bar{e}}$
for each $e \in E(A)$, and injections
$\alpha_e \colon \A_e \inclusion \A_{\boundary_1(e)}$
from each edge group into the group at its terminal vertex.
If $v_0 \in V(A)$ then the \emph{fundamental group}
$\Gamma = \pi_1(\mathbf{A},v_0)$ acts without edge inversions
on the \emph{universal covering tree}
$(\widetilde{\mathbf{A},v_0})$ with quotient~$A$,
so that for each lift $\tilde{e}$ of an edge~$e$
the inclusion $\alpha_e$ is isomorphic to the inclusion
$\Fix_\Gamma(\tilde{e}) \inclusion \Fix_\Gamma(\boundary_1 \tilde{e})$.
Furthermore, every action (without inversions) arises
in this manner from a quotient graph of groups
$\Gamma \quotient X = (\Gamma \backslash X, \A)$.

Let $i(e) = [\A_{\boundary_1(e)} : \alpha_e \A_e]$,
and let $I(\mathbf{A}) = (A,i)$ denote the corresponding
\emph{edge-indexed graph}.
A \emph{grouping} of an edge-indexed graph $(A,i)$ is a graph
of groups~$\mathbf{A}$ with $I(\mathbf{A}) = (A,i)$.
The grouping is \emph{finite} if all the associated groups are finite
and \emph{faithful} if the action on the universal covering tree
is faithful.

Note that the universal covering tree $(\widetilde{\mathbf{A},v_0})$
depends only on the underlying edge indexed graph $I(\mathbf{A})$.
In particular, if $\boundary_1(e)=v$
then for each lift $\tilde{v}$ of~$v$ there are exactly $i(e)$ lifts of~$e$
with terminal vertex~$\tilde{v}$.

In order to construct groupings of edge-indexed graphs, a useful
intermediate step is an ordering.
An \emph{ordering} of an edge-indexed
graph $(A,i)$ is a function $N \colon V(A) \amalg E(A) \to \Q_+$
such that for each $e\in E(A)$ we have
\[
   N(e) = N(\bar{e}) = \frac{N\bigl(\boundary_1(e)\bigr)}{i(e)}
\]
An edge indexed graph is \emph{unimodular} if it admits an ordering.
Any two orderings of $(A,i)$ differ by a constant multiple,
and thus an ordering is 
uniquely determined by its value at a single vertex.  

An integral valued ordering gives a set of numbers that are combinatorially
admissible as the orders of vertex and edge groups of a finite grouping
of $(A,i)$.
More precisely, if $\mathbf{A}=(A,\A)$ is a finite grouping of $(A,i)$
then the orders of the vertex and edge groups define an integral
ordering $N$ with $N(v) = \abs{\A_v}$.

Conversely, any integral ordering~$N$ has
a finite cyclic grouping, constructed as follows.
Let each vertex group [resp.\ edge group] be finite cyclic with order
given by the value of~$N$ on that vertex [resp.\ edge],
and let $\alpha_e \colon \A_e \inclusion \A_{\boundary_1(e)}$ be given by
$[1] \mapsto \bigl[i(e)\bigr]$.
This grouping is effective if and only if the values
of~$N$ do not have a common factor.

\subsection{Haar measure and covolume}
\label{subsec:HaarMeasure}

Let $G$ be a locally compact topological group with a left invariant
Haar measure~$\mu$.
A discrete subgroup $\Gamma \le G$
is a \emph{$G$--lattice} if the covolume
$\mu(\Gamma \backslash G)$ is finite.

Suppose $G$ acts on a set~$S$ with compact open stabilizers $G_s$
for each $s \in S$.
Then every discrete subgroup $\Gamma \le G$ has finite stabilizers 
$\Gamma_s$,
and we define the \emph{$S$--covolume} of $\Gamma$ to be the quantity
\[
   \Vol_S (\Gamma \quotient S)
     := \sum_{s\in \Gamma \quotient S} 1/ \abs{\Gamma_s}
\]

The following theorem
relates Haar measure and $S$--covolume of discrete subgroups.

\begin{theorem}[\cite{BL}, Chapter~1]\label{thm:Haar}
Let $G$ be a locally compact topological group acting on a set~$S$ with
compact open stabilizers and a finite quotient $G \backslash S$.
Suppose further that $G$ admits at least one lattice.
Then there is a normalization of the Haar measure~$\mu$,
depending only on the choice of $G$--set~$S$,
such that for each discrete subgroup $\Gamma \le G$ we have
\[
   \Vol_S(\Gamma \quotient S) = \mu(\Gamma \backslash G)
\]
\end{theorem}

When $G$ acts by automorphisms on a tree~$X$,
the usual convention in the literature is to
compute covolumes of lattices with respect to the action on the set
of all vertices of~$X$.
However, other natural choices for~$S$ often exist.
For instance, one could let $S$ be the set of edges of~$X$,
or $S$ could be a union of $G$--orbits of vertices.

The preceding theorem shows that the particular choice of~$S$
is largely a matter of convenience, as long as the same choice is
made for all computations.

\section{Commensurability invariants}
\label{sec:Invariants}

In this section we discuss and 
construct several commensurability invariants for
discrete subgroups of $\Aut(X)$, for any locally finite tree~$X$.  
The invariants include the covolume, the growth of the quotient graph,
and the growth of the orders of vertex stabilizers in the quotient
graph of groups.
We remark that each of the invariants
in this section can be easily extended to commensurability invariants 
of lattices in the automorphism group of any locally finite simplicial 
complex.  

\subsection{Coverings of graphs of groups}
\label{subsec:Covers}

Suppose $X$ is a locally finite tree and $H$ a discrete subgroup
of $\Aut(X)$ acting without inversions on~$X$.
As explained below,
a subgroup $\Gamma \le H$ corresponds to a covering
$\Gamma \quotient X \to H \quotient X$ of the corresponding
graphs of groups.
In this subsection we examine some features of this correspondence,
focusing on the case when $\Gamma$ has finite index in~$H$.
Related results can be found in \cite{BL} and \cite{Ro}.

The notion of a covering of graphs of groups has been defined by Bass
(\cite{Ba}) and also by Bridson--Haefliger in the context of complexes
of groups (\cite{BH99}).
The exact definitions given by Bass and Bridson--Haefliger
are rather technical and will not be given here.
In fact the two definitions are each sufficiently technical that
it is not currently known whether they are equivalent,
although their equivalence is conjectured by Bridson--Haefliger.
We will not need to use the full technical strength of these definitions.
We only use the following weaker facts,
which are a common subset of the definitions of both Bass and
Bridson--Haefliger and are sufficient for our purposes.

Let $\mathbf{A} = (A,\A)$ and $\mathbf{B} = (B, \B)$ denote the
graphs of groups $H \quotient X$ and $\Gamma \quotient X$, respectively,
and let $(A,i)$ and $(B,j)$ be the underlying edge-indexed graphs.
On the level of graphs, the quotient map $X \to A$ factors through~$B$.
Let $q \colon B \to A$ denote the map arising from this factorization.
Now for each $b\in V(B)$ let $a = q(b)$ and choose a lift
$\tilde{b}$ of~$b$ to $X$.
Then the stabilizer $\Fix_{\Gamma}(\tilde{b})$
is naturally a subgroup of $\Fix_H(\tilde{b})$,
so $q$ induces a monomorphism $\B_b \inclusion \A_a$ of vertex groups,
which is well-defined up to conjugacy.
Similarly, $q$ induces monomorphisms of edge groups.

Furthermore, it is easy to see that the induced map
$q \colon (B,j) \to (A,i)$ is a \emph{covering of edge-index graphs}
in the sense of Bass--Lubotzky (\cite{BL}).
In other words, for each vertex $b \in V(B)$ with $q(b)=a$
and each edge~$e$ with $\boundary_1(e) = a$, we have
\begin{equation}
\label{eqn:IndexedCover}
   i(e) = \sum_{f\in q_b^{-1}(e)} j(f)
\end{equation}
where $q_b$ denotes the restriction of $q$ to the edges in
$\boundary_1^{-1}(b)$.

\begin{defn}
The \emph{degree} of the cover $q\colon \mathbf{B}\to\mathbf{A}$
is the quantity
\begin{equation}
\label{eqn:Degree}
   \deg(q) := \sum_{b\in q^{-1}(a)} [\A_a : \B_b]
\end{equation}
for any fixed vertex
$a \in V(A)$.\endrem
\end{defn}

\begin{proposition}\label{prop:Degree}
The degree of~$q$ does not depend on the choice of $a \in V(A)$.
Furthermore, if $H$ and $\Gamma$ each have finite covolume, then
$\deg(q) = [H:\Gamma]$.
\end{proposition}

For instance, if each inclusion $\B_b\inclusion \A_a$ is an isomorphism,
then each point of~$A$ lifts to exactly $\deg(q)$ points in~$B$,
and the graph map $q\colon B\to A$ is a covering of spaces
whose topological degree is equal to $\deg(q)$.
Alternately, if the graph map $q$ is an isomorphism, then each point
of~$A$ has exactly one lift to~$B$, and we have
$[\A_a : \B_b] = \deg(q)$ whenever $q(b)=a$.

\begin{proof}[Proof of Proposition~\ref{prop:Degree}]
Applying the definition of $i$ and $j$ to (\ref{eqn:IndexedCover}), gives
\[
   \frac{\abs{\A_a}}{\abs{\A_e}}
     = \sum_{f \in q_b^{-1}(e)} \frac{\abs{\B_b}}{\abs{\B_f}}
\]
or equivalently
\[
   [\A_a : \B_b] = \sum_{f\in q_b^{-1}(e)} [\A_e : \B_f]
\]
Summing over all vertices $b$ with $q(b)=a$ gives
\[
   \sum_{b\in q^{-1}(a)} [\A_a : \B_b]
     = \sum_{f \in q^{-1}(e)} [\A_e : \B_f]
\]
Notice that the quantity on the left hand side is independent
of the choice of $e \in \boundary_1^{-1}(a)$,
and the right hand side is unchanged if we replace $e$ with $\bar{e}$.  
Hence both boundary vertices of an edge give the same degree.
Since $A$ is connected,
it follows that $\deg(q)$
is independent of the choice of $a \in V(A)$.

Dividing both sides of (\ref{eqn:Degree}) by $\abs{\A_a}$
gives
\[
   \frac{\deg(q)}{\abs{\A_a}} = \sum_{b\in q^{-1}(a)} \frac{1}{\abs{\B_b}}
\]
which immediately leads to the formula
$\deg(q) \Vol(\mathbf{A}) = \Vol(\mathbf{B})$.
(Here volumes are computed with respect to the set of all vertices.)
Thus
\[
   \deg(q) = \frac{\mu(\Gamma\backslash G)}{\mu(H\backslash G)}
     = [H:\Gamma],
\]
where the second equality is a standard fact about lattices in locally
compact groups (see, \textit{e.g.}, \cite{Rag}).
\end{proof}

\subsection{Quotient growth type}
\label{subsec:QuotientGrowth}

Let $f,g\colon \N \to \N$ be any two functions.
We say that $f\preceq g$ if
\[
   f(k) \le \alpha\,g(k+\beta)
\]
for some $\alpha,\beta\in \N$.
We say $f$ and~$g$ are \emph{equivalent}, denoted $f \simeq g$, 
if $f \preceq g$ and $g \preceq f$.
Note that there is no
multiplicative factor inside the argument of $g$.

Let $A$ be a locally finite graph with basepoint~$\star$, and let 
$g\colon \N \to \N$ be such that $g(k)$ the number of vertices
in the combinatorial ball of radius~$k$ centered at~$\star$.  
We define the \emph{growth type} of $(A,\star)$ to be the 
equivalence class of the function $g$.  
It is easy to see that changing the basepoint does not change the 
growth type.  

Now suppose $\Gamma \le H$
are two discrete subgroups of $\Aut(X)$ with $[H : \Gamma] < \infty$.
Let $q\colon \Gamma\quotient X \to H\quotient X$ be the
corresponding finite degree covering map.
The following result is an immediate corollary of
Proposition~\ref{prop:Degree}.

\begin{corollary}
For each vertex $a \in V(H\backslash X)$, the set
$q^{-1}(a) \subseteq V(\Gamma\backslash X)$ contains at most
$\deg(q)$ vertices.
\qed
\end{corollary}

The following is an easy exercise using the previous corollary.

\begin{proposition}
Let $X$ be a locally finite tree, and
let $\Gamma$ and~$\Gamma'$ be commensurable lattices of $\Aut(X)$.
Then $\Gamma\backslash X$ and $\Gamma'\backslash X$ have the same growth
type.
\qed
\end{proposition}

\subsection{Stabilizer growth type}
\label{subsec:StabilizerGrowth}

Let $\mathbf{A} = (A,\A)$ be a graph of finite groups such that $A$ is
a locally finite graph.
Fix a basepoint $\star \in V(A)$.
The \emph{stabilizer growth type} of ($\mathbf{A},\star)$ is the
equivalence class of the function $g\colon \N \to \N$ such that
$g(k)$ is the order of the largest vertex group associated to any
vertex in the combinatorial ball of radius~$k$ centered at~$\star$.
Clearly the equivalence class of~$g$ does not depend on the choice
of basepoint~$\star$.

A finer invariant of~$\mathbf{A}$ is the \emph{\boldmath$p$--stabilizer growth
type} for any fixed prime~$p$, defined as follows.
If a finite group~$K$ has order $p_1^{n_1}\dotsm p_\ell^{n_\ell}$
for distinct primes $p_1,\dots,p_\ell$,
then the \emph{\boldmath$p_i$--order} of~$K$, denoted $\abs{K}_{p_i}$,
is the natural number $p_i^{n_i}$.  
We define the \emph{\boldmath$p$--stabilizer growth type} of
$(\mathbf{A},\star)$, for $p$ prime, as the equivalence class of the
function $g\colon \N\to\N$ such that $g(k)$ is the maximal $p$--order
of any vertex group associated to a vertex in the combinatorial ball
of radius~$k$ centered at~$\star$.
The $p$--stabilizer growth type also does not depend on the choice of
basepoint.

Suppose $\Gamma \le H$ are two discrete subgroups of $\Aut(X)$,
for some locally finite tree~$X$, such that $[H:\Gamma] < \infty$,
and let $\mathbf{A}=(A,\A)$ and $\mathbf{B}=(B,\B)$ be the graphs
of groups $H\quotient X$ and $\Gamma\quotient X$ respectively.
Let $q\colon \mathbf{B} \to \mathbf{A}$ denote the associated
covering map.
The following lemma is an immediate corollary of
Proposition~\ref{prop:Degree}.

\begin{lemma}
For each $a \in V(A)$, each $b \in q^{-1}(a)$,
and each prime~$p$, the order \textup{[}resp.\ $p$--order\textup{]}
of\/ $\B_b$ divides the order \textup{[}$p$--order\textup{]} of $\A_a$.
Furthermore, we have
\[
   1 \le \frac{\abs{\A_a}}{\abs{\B_b}} \le \deg(q)
     \quad \text{and} \quad
   1 \le \frac{\abs{\A_a}_p}{\abs{\B_b}_p} \le \deg(q).
\]
\qed
\end{lemma}

\begin{corollary}
Let $X$ be a locally finite tree, and let $\Gamma$ and $\Gamma'$
be commensurable lattices in $\Aut(X)$.
Then $\Gamma\quotient X$ and $\Gamma'\quotient X$ have the same
stabilizer growth type and also the same $p$--stabilizer
growth type for each prime~$p$.\qed
\end{corollary}

\begin{remark}
Other related commensurability invariants exist using variations
on the preceding definitions.
For instance, combinatorial spheres are coarsely preserved by finite
covers.
So one could consider the growth of the minimal order, or minimal $p$--order,
of a vertex group in the combinatorial sphere.

More generally, given a discrete $\Gamma \le \Aut(X)$
there is a spectrum of orders, and $p$--orders,
of vertex groups in each combinatorial sphere of~$X$,
which is coarsely preserved under commensurability.
It seems possible that a weighted variant of
Patterson--Sullivan measure on $\boundary X$ could be associated to each 
commensurability class of discrete subgroup $\Gamma\le \Aut(X)$.
The idea is to weight each vertex by the relative size of the
associated vertex group compared to other vertex groups in the same
sphere.
\endrem
\end{remark}

\section{Constructions of lattices}
\label{sec:construct}

As mentioned in the introduction, Bass--Lubotzky \cite{BL} (for $m=n$)
and Rosenberg \cite{Ro} (for every $m\geq n\geq 3$) proved that for
every real number $r>0$, there exists a nonuniform lattice $\Gamma$ in
$G=\Aut(X_{m,n})$ with covolume $r$.
In this section we give a different
construction which will allow us the flexibility of exhibiting the
lattices promised in Theorem \ref{thm:QuotientGrowth}.  

The construction will occur in three steps: the basic building block
will be given in \S\ref{subsec:StarTrees}; these will be used in 
\S\ref{subsec:ArbitraryVolumes} to give 
lattices with arbitrary covolume bigger than a certain constant
(depending on $m,n$); and the constant will be improved to zero in 
\S\ref{subsec:SmallVolumes}.

\subsection{Star trees and their associated lattices}
\label{subsec:StarTrees}

Let $X_{m,n}$ denote the biregular tree with degrees $m$ and~$n$,
for $3\le n \le m$.
In this subsection we construct a family of discrete subgroups of
$G=\Aut(X_{m,n})$ which will be the basic building blocks for our
construction of the lattices.  
Each subgroup~$\Gamma$ will constructed as the fundamental group of a graph
of groups $\mathbf{A} = (A,\A)$ with $A$ a tree.

A \emph{star of degree~\boldmath$m$} is a tree with
one vertex of degree~$m$, called the \emph{center} of the star,
and $m$ vertices of degree one.
A \emph{star tree of degree~\boldmath$m$} is a bipartite tree with vertex sets
$V_0$ and~$V_1$,
such that every vertex of~$V_0$ has degree~$m$ and every vertex of~$V_1$
has degree either one or two.
Any star tree of degree~$m$ can be uniquely expressed as a union of stars of degree~$m$ provided $m \ge 3$.

Let $A$ be a star tree of degree $m$ with basepoint $v_0 \in V_0$.
Put the following edge indexing on~$A$.
For each edge $e$ such that $\boundary_1(e) \in V_0$, set $i(e) = 1$.
For each edge $e$ such that $\boundary_1(e)$ is a vertex of~$V_1$ with
degree one, set $i(e) = n$.
If $v \in V_1$ has degree two, let $e$ and~$e'$ be the two edges
with terminal vertex~$v$, chosen so that $e$ is closer than~$e'$ to~$v_0$ .
Let $i(e) = n-1$ and $i(e')=1$.
Note that the edge-indexed graph $(A,i)$ has universal cover~$X_{m,n}$.

Since $A$ is a tree, any edge indexing on~$A$ is unimodular.
In particular, if we let $N$ be the ordering of $(A,i)$
normalized so that $N(v_0) = 1$,
it is easy to see that $N$ is integral.

The \emph{canonical grouping} of $(A,i)$ is the cyclic grouping~$\mathbf{A}$
associated to~$N$, which is necessarily faithful,
since each lift of $v_0$ has a trivial stabilizer.
The canonical grouping is a graph of finite groups~$\mathbf{A}$
with universal cover $X_{m,n}$.

\begin{convention}[Star tree covolumes]
For convenience, we compute all covolumes of star tree lattices
with respect to the set of vertices~$V_0$, using Theorem~\ref{thm:Haar}.
(In case $m=n$, we first pass to the index two
subgroup $G_0 < G$ stabilizing the lift of~$V_0$.)
For instance, if $A$ consists of just a single star, then its canonical
grouping~$\mathbf{A}$ has $V_0$--covolume $1$.
On the other hand, the covolume with respect to the full set of vertices
$V = V_0 \coprod V_1$ is $(n+m)/n$.
Thus the reader can multiply any $V_0$--covolume by the factor $(n+m)/n$
to obtain the corresponding $V$--covolume.
\end{convention}

For each $v \in V_0$, the \emph{level} of~$v$, denoted $\ell(v)$, 
is the number of vertices of $V_1$ on the unique edge-path from $v$
to the basepoint~$v_0$.
If $\mathbf{A}$ is the canonical grouping~$\mathbf{A}$ of
a star tree~$A$, then for each $v \in V_0$, we have
$N(v) = (n-1)^{\ell(v)}$.
Therefore the covolume of~$\mathbf{A}$ is given by
\[
    \sum_{v\in V_0} \frac{1}{(n-1)^{\ell(v)}}
\]
establishing the following theorem.

\begin{theorem}
\label{thm:VolumeFormula}
Suppose $3\le n\le m$,
and let $A$ be any star tree of degree~$m$.
Then there is a discrete subgroup $\Gamma(A) \le \Aut(X_{m,n})$,
such that $\Gamma(A) \backslash X_{m,n} = A$, and $\Gamma(A)$ has covolume
\[
   \sum_{v\in V_0(A)} \frac{1}{(n-1)^{\ell(v)}}
\]
\end{theorem}

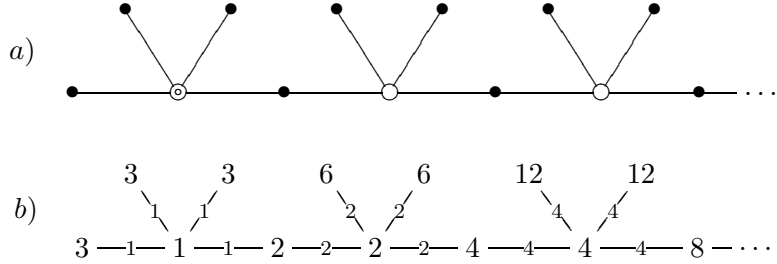
\begin{figure}[t]
\medskip
\[
a) \quad
  \vcenter{
   \xymatrix@R=25pt@C=14pt{
       & *=0{\bullet}\ar@{-}[dr]
       & & *=0{\bullet}\ar@{-}[dl]
       & & *=0{\bullet}\ar@{-}[dr]
       & & *=0{\bullet}\ar@{-}[dl]
       & & *=0{\bullet}\ar@{-}[dr]
       & & *=0{\bullet}\ar@{-}[dl] \\
      *=0{\bullet}\ar@{-}[rr]
      & & *+[o][F=]{}\ar@{-}[rr]
      & & *=0{\bullet}\ar@{-}[rr]
      & & *+[o][F]{}\ar@{-}[rr]
      & & *=0{\bullet}\ar@{-}[rr]
      & & *+[o][F]{}\ar@{-}[rr]
      & & *=0{\bullet}\ar@{-}[r]
      & {\dots}
   }
  }
\]
\bigskip
\[
  b) \quad
  \vcenter{
   \xymatrix@R=15pt@C=7pt{
       & {3}\ar@{-}|{1}[dr]
       & & {3}\ar@{-}|{1}[dl]
       & & {6}\ar@{-}|{2}[dr]
       & & {6}\ar@{-}|{2}[dl]
       & & {12}\ar@{-}|{4}[dr]
       & & {12}\ar@{-}|{4}[dl] \\
      {3}\ar@{-}|{1}[rr]
      & & {1}\ar@{-}|{1}[rr]
      & & {2}\ar@{-}|{2}[rr]
      & & {2}\ar@{-}|{2}[rr]
      & & {4}\ar@{-}|{4}[rr]
      & & {4}\ar@{-}|{4}[rr]
      & & {8}\ar@{-}[r]
      & {\dots}
   }
  }
\]
\caption{(a) A star ray~$R$ of degree~$4$.
The basepoint is indicated by a double circle.
(b) The ordering of the canonical grouping of $(R,i)$ in the case
$m=4$ and $n=3$.}
\label{fig:StarRay}
\end{figure}

\begin{example}
\label{exmp:StarRay}
A \emph{star ray} is a minimal infinite star tree.
It consists of an infinite sequence of stars connected end to end.
A star ray~$R$ of degree~$4$ is illustrated in
Figure~\ref{fig:StarRay}(a).
The ordering of the canonical grouping of $(R,i)$ is shown in
Figure~\ref{fig:StarRay}(b) in the case $m=4$ and $n=3$.
This grouping has covolume
\[
   \sum_{\ell=0}^{\infty} \frac{1}{(n-1)^\ell} = \frac{n-1}{n-2} = 2
\]
\end{example}

\subsection{Arbitrary volumes bounded away from zero}
\label{subsec:ArbitraryVolumes}

In this subsection we construct uncountably many lattices in~$G$
with covolume~$\kappa$, for each $\kappa>\kappa_0 := (n-1)/(n-2)$.

If the canonical grouping of a star tree~$A$ is to be a nonuniform
lattice, it is necessary that $A$ contains at least one ray.
Thus the
covolume of any nonuniform lattice produced in this manner is bounded
below by $\kappa_0$, which is the covolume associated to a
star ray~$R$ (see Example~\ref{exmp:StarRay}).

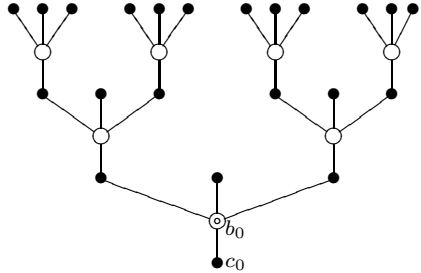
\begin{figure}[t]
\[
  \xymatrix@R=10pt@C=5pt{
     *=0{\bullet}\ar@{-}[dr] & *=0{\bullet}\ar@{-}[d]
       & *=0{\bullet}\ar@{-}[dl]
     & & *=0{\bullet}\ar@{-}[dr] & *=0{\bullet}\ar@{-}[d]
       & *=0{\bullet}\ar@{-}[dl]
     & & *=0{\bullet}\ar@{-}[dr] & *=0{\bullet}\ar@{-}[d]
       & *=0{\bullet}\ar@{-}[dl]
     & & *=0{\bullet}\ar@{-}[dr] & *=0{\bullet}\ar@{-}[d]
       & *=0{\bullet}\ar@{-}[dl]                               \\
  &  *+[o]
  [F]{}\ar@{-}[d] & & & & *+[o][F]{}\ar@{-}[d]
    & & & &  *+[o][F]{}\ar@{-}[d] & & & & *+[o][F]{}\ar@{-}[d]   \\
  & *=0{\bullet}\ar@{-}[drr] & & *=0{\bullet}\ar@{-}[d] & &
    *=0{\bullet}\ar@{-}[dll] & &
    & & *=0{\bullet}\ar@{-}[drr] & & *=0{\bullet}\ar@{-}[d] & &
    *=0{\bullet}\ar@{-}[dll]                                     \\
  & & & *+[o][F]{}\ar@{-}[d] & & & & & & & & *+[o][F]{}\ar@{-}[d] \\
  & & & *=0{\bullet}\ar@{-}[drrrr] & & & & *=0{\bullet}\ar@{-}[d]
     & & & & *=0{\bullet}\ar@{-}[dllll]                          \\
  & & & & & & & *+[o][F=]{}\ar@{-}[d]^>{c_0}^<{b_0}              \\
  & & & & & & & *=0{\bullet}
    }
\]
\caption{The tree $B_3$ in the case $m=4$ and $n=3$.}
\label{fig:FullTree}
\end{figure}

For each natural number~$p$, let $B_p$ be a finite star tree of degree~$m$
with basepoint $b_0$
such that $V_0(B_p)$ contains exactly $(n-1)^j$ vertices of level~$j$
for $j = 0,\dots,p-1$ and no vertices of any higher level.
Let $c_0 \in V_1(B_p)$ be a vertex of degree one adjacent to~$b_0$.
The tree $B_3$ is illustrated in Figure~\ref{fig:FullTree}
in the case $m=4$, and $n=3$.

Consider the star tree obtained by gluing a copy of $B_p$ to the star
ray~$R$ along the vertex $c_0$ so that $b_0$ has level~$j$ with respect
to the basepoint of~$R$.
We refer to this process as \emph{gluing $B_p$ to~$R$ at level~$j$}.
Such a gluing increases the covolume by precisely
\[
   \frac{p}{(n-1)^j}
\]

It is now easy to construct lattices with any given 
covolume $\kappa > \kappa_0$, as follows.
Choose a sequence of natural numbers $(e_j)$ such that
\[
  \rho := \kappa-\kappa_0
   = \frac{e_1}{n-1} + \frac{e_2}{(n-1)^2} + \dots + \frac{e_j}{(n-1)^j}
       + \dotsb
\]
and let $A$ be the tree obtained by gluing for each $j>0$ a copy
of $B_{e_j}$ to~$R$ at level~$j$.
The canonical grouping of $A$
is then a lattice $\Gamma$ with covolume~$\kappa$.

Choosing the sequence $(e_j)$ is similar to writing $\rho$ in base $n-1$,
except we do not require that each $e_j$ be smaller than $(n-1)$.
Consequently, there is a large amount of flexibility in choosing the~$e_j$,
even if we restrict ourselves to bounded sequences.
In particular, by varying the sequence $(e_j)$, it is easy to obtain
uncountably many lattices with a given covolume.
We will take advantage of this flexibility below

\subsection{Arbitrarily small volumes}
\label{subsec:SmallVolumes}

In this subsection we modify the lattices constructed above to produce
lattices with arbitrarily small covolumes.

Recall that a star tree~$A$ of degree~$m$ gives rise to an edge-indexed
graph $(A,i)$ with universal cover $X_{m,n}$.
Thus far, we have only considered the covolume of the canonical grouping
of $(A,i)$.
Consequently, we have not produced lattices with covolume
arbitrarily close to zero.
In order to achieve smaller covolumes, we must increase the orders of
the edge and vertex groups.

The main source of tension in increasing the orders of the local groups
is that the grouping must remain faithful.
If a grouping $(\mathbf{A},\A)$
is not faithful, then the kernel of the action on the universal cover
is a group~$K$ which can be embedded as a nontrivial normal subgroup in
every edge group and vertex group
in a manner which is equivariant with respect to the edge maps
(see \cite[I.1.23]{Ba}).

Loosely speaking, a nonuniform lattice arising as the canonical grouping
of a star tree is based on a sequence of inclusions
\[
  \xymatrix{
     G_0 \ar[r]^{\iota_0} & G_1 \ar[r]^{\iota_1}
     & {\dotsb} \ar[r]^{\iota_{j-1}}
     & G_j \ar[r]^{\iota_j} & {\dotsb}
  }
\]
where $\abs{G_j} = (n-1)^j$.
Our strategy now is to generalize the canonical grouping by
finding a group $H$ and a sequence of maps
$\phi_j \colon H \to \Aut(G_j)$ equivariant with respect to the~$\iota_j$.
The $\iota_j$ then induce a sequence of inclusions
\[
  \xymatrix{
     G_0 \sdirect_{\phi_0} H \ar[r] & G_1 \sdirect_{\phi_1} H \ar[r]
     & {\dotsb}\ar[r]
     & G_j \sdirect_{\phi_j} H \ar[r] & {\dotsb}
  }
\]
Since $G_0$ is the trivial group,
the corresponding grouping is faithful if and only if $H$ does not have 
a subgroup~$H'$ which is normal in every semidirect product,
or equivalently if for each nontrivial element $h \in H$
there is an $i$ such that $h$ acts nontrivially on $G_i$.

Note that the volume of this new lattice equals $\rho/\abs{H}$,
where $\rho$ is the volume of the original lattice, since the
new ordering differs from the old only by the multiplicative factor
$\abs{H}$.
To get arbitrarily small covolumes, we want to be able to
choose the group~$H$ arbitrarily large.  

It remains to produce specific groups $H$ and $G_j$
and actions $H \to \Aut(G_j)$ satisfying the preceding properties.
Toward this end, fix a natural number $k$, and define
\[
   G_j = \begin{cases}
           \Z/(n-1)^j\Z, & \text{if $j \le k$} \\
	   \Z/(n-1)^k\Z \times \Z/(n-1)^{j-k}\Z, & \text{if $j > k$}
         \end{cases}
\]
Let $\iota_j \colon G_j \to G_{j+1}$ be defined as follows.
If $j<k$, then $\iota_j$ is multiplication by $n-1$, and if
$j\ge k$, it is the identity on the first factor and multiplication by
$n-1$ on the second factor.

Now let $H$ be the multiplicative group of units of the ring
$\Z/{(n-1)^k\Z}$.
There is a natural monomorphism
$H \inclusion \Aut(G_k)$ given by multiplication,
which pulls back to maps $H \to \Aut(G_j)$ for each $j\le k$.
For $j\ge k$, let $H$ act by $\phi_k$ on the left factor and trivially
on the right factor.
It is easy to see that this defines
a sequence of maps $\phi_j \colon H \to \Aut(G_j)$ equivariant with
respect to the~$\iota_j$.
Furthermore, by construction,
every element of~$H$ acts nontrivially on $G_k$.
But the order of $H$ can be made arbitrarily large by choosing $k$
sufficiently large.

For each star tree~$A$ of degree~$m$,
the preceding discussion can be applied to produce finite
groupings of $(A,i)$ with arbitrarily small covolumes as follows.
For each vertex $v \in V_0$, set $\A_v = G_{\ell(v)} \sdirect H$.
If $e$ is an edge with $\boundary_1(e) = v$, set $\A_e = \A_v$.
If $w \in V_1$ has degree one and $e$ is the edge with $\boundary_1(e) = w$, set $\A_w = \A_e \times \Z/n\Z$ and let the edge map be the inclusion
of the first factor.
If $w \in V_1$ has degree two, then $\boundary_1^{-1}(w)$ consists
of two edges $e_1$ and $e_2$ such that $\A_{e_1} = G_i \sdirect H$
and $\A_{e_2} = G_{i+1} \sdirect H$.
In this case,
set $\A_w = \A_{e_2}$, and let the map $\A_{e_1} \inclusion \A_{w}$
be the map $G_i \sdirect H \inclusion G_{i+1} \sdirect H$ induced by
$\iota_i$ as above.

By the argument above, this construction produces, for each~$k$,
a finite grouping $\mathbf{A}_k$ of $(A,i)$.
The grouping $\mathbf{A}_k$ is faithful provided that the tree~$A$
contains at least one vertex at level~$k$.
Thus if $A$ is infinite, $\mathbf{A}_k$ is faithful for every~$k$.

\section{Realizing different values of the invariants}
\label{sec:construct2}

In this section we state the precise version of Theorem
\ref{theorem:informal} (see Theorem \ref{theorem:newmain} below) and 
use the constructions given in \S\ref{sec:construct} to exhibit its proof.

\subsection{Different growth types and covolumes}
\label{sec:DifferentGrowth}

Throughout this section, all
quotient growths should be computed with respect to
the metric on $X_{m,n}$ in which each edge has length $1/2$.
In this metric, the function $\ell$ measures the distance from a vertex
to the basepoint.

A function $f\colon \N \to \N$ is an \emph{acceptable quotient growth
function} if it satisfies
\[
   f(0)=1, \quad 1 \le f(j+1) \le 2 f(j), \quad
   \text{and} \quad
   \sum_{j=0}^\infty \frac{f(j)}{2^j} < \infty
\]
(The third condition is a nontrivial requirement only in the case when $n=3$.)

\begin{lemma}\label{lem:Acceptible}
Let $f$ be any acceptable quotient growth function.
Then there is a star tree $T_f$ of degree~$m$
whose growth is equivalent to~$f$.
Furthermore, we may assume that $T_f$ has a vertex $c_0 \in V_1(T_f)$
of degree one adjacent to the basepoint.
\end{lemma}

\begin{proof}
Construct $T_f$ with exactly $f(j)$ stars at level $j$.
The inequality $1 \le f(j+1) \le 2f(j)$ quarantees that this can always
be done, since each level contains at most twice the number of stars
as the previous level and $2 \le m-1$.
\end{proof}

\begin{theorem}\label{thm:QuotientGrowth}
Let $f$ be any acceptable quotient growth function,
and $\kappa$ any positive real number.
Then there is a lattice $\Gamma$ in $\Aut(X_{m,n})$
with covolume~$\kappa$ and quotient growth type~$f$.
In particular, there are uncountably many commensurability classes
of lattices with covolume~$\kappa$, and also uncountably many
commensurability classes with growth type~$f$.
\end{theorem}

\begin{proof}
We will construct our lattices as canonical groupings of star trees.
Without loss of generality, we may assume that $\kappa >\kappa_0$,
since the techniques of \S\ref{subsec:SmallVolumes}
can then be used to obtain arbitrary positive covolumes.

Let $T_f$ be a star tree of degree~$m$ with growth equivalent
to~$f$, as given by Lemma~\ref{lem:Acceptible}.
The canonical grouping of such a tree is a lattice with finite
covolume~$\nu$.

Fix a natural number $k$ sufficiently large that
$\nu /(n-1)^k < \kappa - \kappa_0$,
and choose a bounded sequence of natural numbers $(e_j)_{j \ne k}$ so that
\[
  \kappa - \kappa_0
          = \frac{e_1}{n-1} + \dots
          + \frac{e_{k-1}}{(n-1)^{k-1}}
	  + \frac{\nu}{(n-1)^k}
          + \frac{e_{k+1}}{(n+1)^{k+1}} + \dotsb
\]

Form a star tree $A(f,\kappa)$ by gluing $T_f$ to the star ray~$R$
at level~$k$
and for each $j\ne k$ gluing a copy of $B_{e_j}$ at level~$j$.
Since the $B_{e_j}$ have uniformly bounded depth, they do not affect the
growth type of~$A$, which is therefore equivalent to~$f$.
Furthermore, the canonical grouping of~$A$ has covolume~$\kappa$
by construction.

The last assertion of the theorem
follows from the fact that for each $\alpha \in (1,2)$
there is an acceptable quotient growth function
equivalent to $j\mapsto \alpha^j$ and two such functions are inequivalent
whenever $\alpha \ne \alpha'$.
\end{proof}

\subsection{Different stabilizer growths}
\label{sec:DifferentStabGrowth}

In the previous sections we constructed lattices in $\Aut(X_{m,n})$
as groupings of the
canonical edge indexing of a star tree.
All nonuniform lattices arising from such a canonical edge indexing
necessarily have the same stabilizer growth type, namely the type of the
exponential function given by $h(k) = (n-1)^{k}$.
In this subsection we construct lattices with different stabilizer growths
by varying the edge indexing.
The techniques of this subsection require that $n$ be composite.

As in the previous subsection, we continue to compute all growth functions
using the metric in which each edge has length $1/2$.

We say that a sequence $s = (s_k)$ is 
\emph{\boldmath$n$--admissible} if each $s_k$ 
is an integer dividing $n$ and if $2<s_k\leq n$.  
For each~$k$, let $r_k := n/s_k$.
Given an $n$--admissible sequence~$s$
and an infinite star tree $A$ of order~$m$ with basepoint $v_0$, 
there is an edge indexing $I(A,s)$ with universal cover $X_{m,n}$
defined as follows.
If $\boundary_1(e) \in V_0(A)$, set $i(e) := 1$.
If $\boundary_1(e) \in V_1(A)$ has degree one, set
$i(e) := n$.
If $w\in V_1$ has degree two, then $\boundary_1^{-1}(w) = \{e,e'\}$,
where $e$ is closer to the basepoint $v_0$ than~$e'$.
Then the vertices $\boundary_0(e)$ and $\boundary_0(e')$
adjacent to~$w$ satisfy
$\ell \bigl(\boundary_0(e) \bigr) = k-1$
and $\ell \bigl(\boundary_0(e') \bigr) = k$ for some $k$.
In this case, set
\[
   i(e) := n - r_k = r_k(s_k-1) \quad \text{and} \quad
   i(e') := r_k
\]
Observe that the canonical edge indexing is the special case
when $s_k = n$ and $r_k = 1$ for all $k$.

Let $N$ be the ordering associated to $I$ such that $N(v_0) = 1$.
If $v,v' \in V_0(A)$
with $\ell(v) = k-1$ and $\ell(v') = k$, then we have
\[
   N(v') = N(v)\,(s_k-1) = (s_1 - 1) \dotsm (s_k - 1)
\]
In particular, $N$ is integral and $N(v)$ depends only on the level
of $v\in V_0(A)$.
The value $N(w)$ for any $w \in V_1$ is within a factor of~$n$ of $N(v)$
for some $v \in V_0$ at a combinatorial distance~$1$ from~$w$.
Thus we may ignore vertices in~$V_1$ when considering the stabilizer
growth type of~$N$.
It is now clear that $N$ has stabilizer growth type equivalent to the
function $h$ given by
\[
   h(k) = (s_1 - 1) \dotsm (s_k - 1)
\]

In particular, if $t>2$ is any proper, nontrivial factor of~$n$ (\textit{i.e.},
not equal to either $1$ or $n$), then for each real number
$\lambda \in [t-1,n-1]$ there is a sequence
$s = (s_k)$ with $s_k \in \{t,n\}$
such that the ordering~$N$ associated to~$s$
has stabilizer growth equivalent to
the function $k \mapsto \lambda^k$.
Furthermore, distinct choices of $\lambda$
give rise to inequivalent growth types.
Thus we have established the following result.

\begin{prop}
Suppose that $3 \le m \le n$ and that $n>4$ is composite.
Then $G=\Aut(X_{m,n})$
admits uncountably many nonuniform lattices with  distinct
stabilizer growth types.
\end{prop}

Notice that the stabilizer growth type of any ordering $N$ associated to an
$n$--admissible sequence~$s$ depends only on the choice of~$s$.
The growth type is independent of the choice of ordering~$N$ of~$I$,
and is also independent of the choice of infinite star tree~$A$.

Our goal for the remainder of this subsection is to modify the constructions
in the previous subsections to produce lattices modeled on an 
arbitrary $n$--admissible sequence~$s$ that also have arbitrary
acceptable quotient growth and arbitrary covolume.
Note that these constructions involve modifying only the tree~$A$
and the grouping of~$I$.
Thus they do not affect the stabilizer
growth of the corresponding lattices.

\begin{theorem}
\label{theorem:newmain}
Suppose that $3 \le m \le n$ and that $n>4$ is composite.
Let $f$ be any acceptable quotient growth function,
$s=(s_k)$ any $n$--admissible
sequence, and $\kappa$ any positive real number.
Then $G = \Aut(X_{m,n})$ contains nonuniform lattices
with quotient growth type~$f$, covolume~$\kappa$, and stabilizer growth type
given by the function
\[
   h(0)=1; \quad h(j) = (s_1 - 1) \dotsm (s_j - 1)
\]
\end{theorem}

\begin{proof}
Let $R$ be the star ray of degree~$m$ with basepoint $v_0$,
let $N$ be the ordering associated to the edge indexing $I(R,s)$
such that $N(v_0) = 1$,
and let
\[
   \kappa_0 := \sum_{j=0}^{\infty} \frac{1}{h(j)}
\]
which is the covolume of the ordering~$N$.

We first consider the special case when $\kappa > \kappa_0$,
as in \S\ref{subsec:ArbitraryVolumes}.
Let $T_f$ be a star tree with growth~$f$, which has a vertex
$c_0 \in V_1(T_f)$ of degree one adjacent to the basepoint.
Gluing $T_f$ onto the star ray~$R$ at any level~$k$,
has the effect of increasing the covolume of the ordering~$N$
by
\[
   \nu_k := \sum_{j=0}^\infty \frac{f(j)}{h(j+k)}
\]
Since $s$ is $n$--admissible, it follows that $s_j-1 \ge 2$, so that
$h(j) \ge 2^j$.
Therefore
\[
   \nu_k \le 2^{-k} \sum_{j=0}^\infty \frac{f(j)}{2^j}
\]
The preceding sum converges since $f$ is an acceptable growth function.
As in the proof of Theorem~\ref{thm:QuotientGrowth},
we fix a natural number~$k$ sufficiently large that
$\nu_k < \kappa - \kappa_0$, and glue $T_f$ to~$R$ at level~$k$.

The next step is to modify the construction of the bounded trees $B_p$
from \S\ref{subsec:ArbitraryVolumes} to take into account the variable
stabilizer growth provided by the sequence~$s$.
For each $p,q \in \N$ let $B_{p,q}$ be a finite depth star tree
of degree~$m$ with exactly $h(q+j)/h(q)$ vertices of level~$j$ for
$j = 0,\dots,p-1$ and no vertices of any higher level.
Thus the process of gluing $B_{p,q}$ to~$R$ at level~$q$
contributes $h(q+j)/h(q)$ vertices of level $q+j$ for each
$j=0,\dots,p-i$.
So $B_{p,q}$ increases the covolume of~$N$ by
\[
   \sum_{j=0}^{p-1} \frac{1}{h(q)} = \frac{p}{h(q)} 
\]

We now construct a star tree $A$ with growth type~$f$ and covolume~$\kappa$
as follows.
Choose a bounded sequence $(e_j)_{j\ne k}$ satisfying
\[
   \kappa - \kappa_0
     = \frac{e_1}{h(1)} + \dots
       + \frac{e_{k-1}}{h(k-1)}
       + \nu_k
       + \frac{e_{k+1}}{h(k+1)} + \dotsb
\]
Then the star tree~$A$ is formed by gluing $T_f$ onto the star ray~$R$
at level~$k$, and for each $j\ne k$ gluing a copy of $B_{e_j,j}$ 
at level~$j$.
The resulting tree~$A$, with the induced ordering~$N$,
has quotient growth type~$f$, covolume~$\kappa$, and stabilizer growth
type~$h$, as desired.

Finally we must consider the case when $\kappa \le \kappa_0$.
But this is easily dealt with by a construction similar to that in
\S\ref{subsec:SmallVolumes}, replacing $(n-1)^j$ with $h(j)$
throughout.
Therefore, given any star tree~$A$ and any $n$--admissible sequence~$s$,
the associated edge indexing $I(A,s)$ admits faithful finite groupings
with arbitrarily small covolumes.
Since quotient growth and stabilizer growth types depend only 
on the underlying edge indexed graph, this process completes the proof
of the theorem.
\end{proof}


\bigskip
\noindent
Department of Mathematics\\
University of Chicago\\
5734 S.~University Ave\\
Chicago, IL 60637, USA\\
E-mails: \texttt{farb@math.uchicago.edu},
\texttt{chruska@math.uchicago.edu}

\end{document}